# Local Rigidity for Simultaneously Diophantine Translations on Tori of Arbitrary Dimension

Boris Petković

**Abstract**

We show that a smooth sufficiently small perturbation of a $\mathbb{Z}^m$ action on the torus by simultaneously Diophantine translations, is smoothly conjugate to the unperturbed action under a natural condition on the rotation sets. This generalizes recent result of Karaliolios [4] of the action generators to higher rank abelian actions, and the result of Moser [8] to higher dimensional tori.

## 1 Introduction

Moser proved in [8], using the rapidly convergent Nash-Moser iteration scheme, that if the rotation numbers of the orientation-preserving commuting diffeomorphisms of the unit circle satisfy a simultaneous Diophantine condition (which we define in section 3) and if the diffeomorphisms are in a $C^\infty$ neighborhood of the corresponding rotations (the neighborhood being imposed by the constants appearing in the simultaneoulsy Diophantine condition), then they are simultaneously smoothly conjugate to rotations. He also showed in [8] that the problem cannot be reduced to that of a single diffeomorphism with a Diophantine rotation number. Namely, there exist an uncountable set of numbers $\gamma_1, \ldots, \gamma_m$ that are simultaneously Diophantine but for all linearly independent vectors $a, b \in \mathbb{Z}^{m+1}$, the ratios

$$\frac{a_0 + a_1\gamma_1 + \ldots + a_m\gamma_m}{b_0 + b_1\gamma_1 + \ldots + b_m\gamma_m}$$

are Liouville numbers. Fayad and Khanin proved in [FK] the same result but in the global setting, i.e. in the case when the commuting circle diffeomorphisms are not necessarily in a small $C^\infty$ neighborhood of the corresponding rotations.

We prove a result similar to that of Moser in [8], but for commuting diffeomorphisms isotopic to the identity on the torus of arbitrary dimension. The lack of a complete invariant for a torus diffeomorphism in dimensions higher then one, such as the rotation number in the case of a circle diffeomorphism, produces difficulties in the study of local and, of course, global rigidity of systems of this type. Instead of a rotation number in dimension one, for a diffeomorphism of a torus there exist in literature at least three different, but related sets generalizing the rotation number in the case of a circle diffeomorphism. They were initially defined by Herman, Misiurewicz and Ziemian in [3] and [6], respectively.

Let $\mathbb{T}^d = \mathbb{R}^d/\mathbb{Z}^d$ be the $d$-dimensional torus, $d \in \mathbb{N}$, and let $\pi\colon \mathbb{R}^d \to \mathbb{T}^d$ be the standard projection map. We use the $\ell^1$ norm on $\mathbb{Z}^d$, i.e., for a vector $k = (k_1, \ldots, k_d) \in \mathbb{Z}^d$ we define its $\ell^1$ norm as $|k| = k_1 + \ldots + k_d$. Let $\mathrm{C}^\infty(\mathbb{T}^d)$ be the space of smooth maps from $\mathbb{T}^d$ to itself and denote by $\mathrm{Diff}_0^\infty(\mathbb{T}^d) \subset \mathrm{C}^\infty(\mathbb{T}^d)$ the set of all diffeomorphisms that are isotopic to the identity map. Now we define $D_j = \frac{\partial}{\partial x_j}$, for $j = 1, \ldots, d$ and $D^k = D_{k_1} \cdots D_{k_d}$, for $k \in \mathbb{N}_0^d$. On $\mathrm{C}^\infty(\mathbb{T}^d)$, and hence on $\mathrm{Diff}_0^\infty(\mathbb{T}^d)$, we introduce a sequence of $\mathrm{C}^r$ norms

$$\|g\|_r = \max_{1 \le j \le d} \max_{|k| \le r} \sup_{x \in \mathbb{T}^d} \left| D^k g_j(x) \right|,$$

where $g = (g_1, \ldots, g_d) \in \mathrm{C}^\infty(\mathbb{T}^d)$ and the maximum in the middle is taken over all vectors $k \in \mathbb{N}_0^d$, for which $|k| \le r$.

As it is very well known, any map $g \in \mathrm{C}^\infty(\mathbb{T}^d)$ can be lifted to a smooth map $G\colon \mathbb{R}^d \to \mathbb{R}^d$, i.e. to a map which satisfies $\pi \circ G = g \circ \pi$. For such a map we have $G(x + k) = G(x) + k$, for all $k \in \mathbb{Z}^d$ and all $x \in \mathbb{R}^d$. Any two such lifts differ by an additive vector in $\mathbb{Z}^d$. For a point $x \in \mathbb{R}^d$, let $\rho(g, x)$ be the set of all limit points of the sequence

$$\left( \frac{G^n(x) - x}{n} \right)_{n=1}^\infty.$$

Since any two lifts of $g$ differ by an integer vector, $\rho(g, x)$ is well defined modulo $\mathbb{Z}^d$. We call

$$\rho_P(g) = \bigcup_{x \in \mathbb{R}^d} \rho(g, x)$$

the pointwise rotation set for $g$. However, there is no reason to fix a point in the last sequence and we now take all the limits of sequences

$$\frac{G^{n_j}(x_j) - x_j}{n_j}, \text{ where } x_j \in \mathbb{R}^d \text{ and } n_j \to \infty, \text{ as } j \to \infty,$$

modulo $\mathbb{Z}^d$, provided they exist. The set obtained in such a way we call the rotation set of $g$ and denote it by $\rho(g)$.

Now we define the third set related to $\rho_p(f)$ and $\rho(f)$. Define a continuous map $\varphi\colon \mathbb{T}^d \to \mathbb{R}^d$ by $\varphi(y) = G(x) - x$, where $x \in \pi^{-1}(y)$. $\varphi(y)$ does not depend on the choice of $x \in \pi^{-1}(y)$. Denote by $\mathcal{M}(g)$ the space of all $g$-invariant probability measures on $\mathbb{T}^d$ and by $\mathcal{M}_e(g)$ its subset consisting of ergodic measures. If $\mu \in \mathcal{M}_e(g)$ then, by Birkhoff ergodic theorem

$$\lim_{n \to \infty} \frac{G^n(x) - x}{n} = \lim_{n \to \infty} \frac{1}{n} \sum_{j=0}^{n-1} \varphi(f^j(y)) = \int_{\mathbb{T}^d} \varphi \, d\mu,$$

for $\mu$-almost every $y \in \mathbb{T}^d$. Denote by $\rho_e(g) = \{\int_{\mathbb{T}^d} \varphi \, d\mu : \mu \in \mathcal{M}_e(g)\}$.

Now that we have defined three different sets that all generalize the rotation number in dimension one, one can ask about the properties of these sets and relations between them. Specifically, one can ask if at least some of them retain some of the properties of the rotation number, for example, the invariance under the conjugation by smooth diffeomorphisms. Obviously, we have inclusions $\rho_e(g) \subset \rho_p(g) \subset \rho(g)$. Misiurewicz and Ziemian proved in [6] that $\rho(g)$ is compact and connected and that $\text{conv}(\rho_e(g)) = \text{conv}(\rho_p(g)) = \text{conv}(\rho(g))$, where $\text{conv}(S)$ denotes the convex hull of a set $S$. Note that $\text{conv}(\rho_e(g)) = \{\int_{\mathbb{T}^d} \varphi \, d\mu : \mu \in \mathcal{M}(g)\}$, since ergodic measures are exactly the extreme points of convex set $\mathcal{M}(g)$. Since we only deal with convex hulls of the rotation sets, we do not prioritize any of these sets as the best way of generalizing the rotation number. Unfortunately, the rotation sets defined above are not, in general, invariant under smooth conjugation. But the rotation set $\rho(g)$ is preserved under the smooth conjugacies isotopic to the identity. This is the reason why we consider only diffeomorphisms isotopic to the identity, instead of, as it is the case in dimension one, the orientation preserving ones. It is also proved in [7] that the rotation set $\rho(g)$ is convex in the two dimensional case.

We use the following convention for constants appearing in estimates that we either use or prove. We denote a positive constant by $C$ and its value can vary from one occurrence to the next. Further, if we want to indicate dependence of a constant on parameters we use subscripts. For example $C_{r,s}$ stands for a positive constant depending on parameters $r, s$. We will not indicate a dependence of a constant on dimension of the torus and on the number of commuting diffeomorphisms, for example, we simply write $C_r$ instead of $C_{d,m,r}$.

We recall the notion of a Diophantine vector. We say that $\alpha \in \mathbb{R}^d$ is Diophantine of type $(\gamma, \tau)$, where $\gamma > 0$ and $\tau > 0$ are constants, if

$$|\langle \alpha, k \rangle - l| > \frac{\gamma}{|k|^\tau}, \tag{1.1}$$

for all $(k, l) \in \mathbb{Z}^d \times \mathbb{Z} \setminus \{(0, 0)\}$.

For any $N \geq 0$, as smoothing operators on $\mathrm{C}^\infty(\mathbb{T}^d, \mathbb{R})$ we use the truncations of Fourier series and denote them by $T_N$. If $g \in \mathrm{C}^\infty(\mathbb{T}^d, \mathbb{R})$ is given by $g(x) = \sum_{k \in \mathbb{Z}^d} \hat{g}(k) e^{2\pi i \langle k, x \rangle}$, then

$$T_N g(x) = \sum_{|k| \leq N} \hat{g}(k) e^{2\pi i \langle k, x \rangle}.$$

Note that $g$ and $T_N g$ have the same average with respect to the Lebesgue measure. We denote it by $Av(g)$. By some abuse of notation, we also use $T_N$ to denote the truncation on $\mathrm{C}^\infty(\mathbb{T}^d, \mathbb{R}^d)$ acting componentwise, and similarly for the average. The following estimates for $g \in \mathrm{C}^\infty(\mathbb{T}^d, \mathbb{R}^d)$ are very well known (see for example [9]).

$$\begin{aligned} \|T_N g\|_r &\leq C_{r,s} N^{r-s+\frac{d}{2}} \|g\|_s \text{ for } 0 \leq s \leq r, \\ \|(I - T_N)g\|_r &\leq C_{r,s} N^{r-s+d} \|g\|_s, \text{ for } 0 \leq r \leq s. \end{aligned} \tag{1.2}$$

## 2 The result

### 2.1 Statement of the Theorem

We discuss the following local question. Let $f_1, \ldots, f_m \in \text{Diff}_0^\infty(\mathbb{T}^d)$ be a finite collection of commuting diffeomorphisms in a $C^\infty$ neighborhood of torus rotations $R_{\alpha_j}(x) = x + \alpha_j \pmod{\mathbb{Z}^d}$, $j = 1, \ldots, m$. Which conditions on vectors $\alpha_1, \ldots, \alpha_m \in \mathbb{R}^d$ imply that there exists $h \in \text{Diff}^\infty(\mathbb{T}^d)$ such that
$$h^{-1} \circ f_j \circ h = R_{\alpha_j}, \text{ for } 1 \leq j \leq m? \qquad (2.1)$$
In this case, the condition for vectors $\alpha_1, \ldots, \alpha_m$ corresponding to the Diophantine condition (1.1), for a single diffeomorphism, is that there exist constants $\gamma > 0$ and $\tau > 0$ such that
$$\max_{j=1\ldots,m} |\langle \alpha_j, k \rangle - l| > \frac{\gamma}{|k|^\tau}, \qquad (2.2)$$
for all $(k,l) \in \mathbb{Z}^d \times \mathbb{Z} \setminus \{(0,0)\}$. If vectors $\alpha_1, \ldots, \alpha_m$ satisfy (2.2) we say that they are simultaneously Diophantine of type $(\gamma, \tau)$. Note that (2.2) is obviously satisfied if at least one of $\alpha'_j s$ satisfies (1.1), for example, $\alpha_1$. Then, the other $\alpha_l$, $l \geq 2$, are completely arbitrary. As we already said, Moser proved in [8] that condition (2.2), in one dimensional case, is strictly weaker then the Diophantine condition (1.1).

We prove the following theorem.

**Theorem 2.1** *Let $\gamma > 0$ and $\tau > d$. There exist $\epsilon = \epsilon(d, \gamma, \tau) > 0$ and $l = l(d, \gamma, \tau) > 0$ such that, if the collection of diffeomorphisms $f_1, \ldots, f_m \in \text{Diff}_0^\infty(\mathbb{T}^d)$ and vectors $\alpha_1, \ldots, \alpha_m \in \mathbb{R}^d$ satisfy*

*(a) $f_j \circ f_l = f_l \circ f_j$, for $1 \leq j, l \leq m$,*

*(b) $\alpha_1, \ldots, \alpha_m$ are simultaneously Diophantine of type $(\gamma, \tau)$,*

*(c) $\max\limits_{j=1,\ldots,m} \left\| f_j - R_{\alpha_j} \right\|_l < \epsilon$*

*(d) $\alpha_j \in \text{conv}(\rho(f_j))$, for $1 \leq j \leq m$,*

*then there exists $h \in \text{Diff}^\infty(\mathbb{T}^d)$ such that*
$$h^{-1} \circ f_j \circ h = R_{\alpha_j}, \text{ for } 1 \leq j \leq m.$$
*Moreover, $h$ can be chosen close to the identity map.*

When $d = 1$, the theorem reduces to the case of commuting circle diffeomorphisms, since in that case the rotation set is just a singleton containing the rotation number. This is also a generalization of Theorem 3.1 in [4]. There, the same theorem is proved but in the case of a single toral diffeomorphism isotopic to the identity map.

## 3 Proof of the result

### 3.1 Notation and conventions

In order to simplify the notation as much as possible, we introduce the following tuples of maps and define their composition, $C^r$ norms, smoothing operators and average componentwise. Namely, if $\mathcal{A} = (a_1, \ldots, a_m)$ and $\mathcal{B} = (b_1, \ldots, b_m)$ are such two tuples, then we define their composition componentwise, i.e.,

$$\mathcal{A} \circ \mathcal{B} = (a_1 \circ b_1, \ldots, a_m \circ b_m),$$

$C^r$ norms as $\|\mathcal{A}\|_r = \max_{1 \leq j \leq m} \|a_j\|_r$, smoothing operators as $T_N \mathcal{A} = (T_N a_1, \ldots, T_N a_m)$ and average as $Av(\mathcal{A}) = (Av(a_1), \ldots, Av(a_m))$. Note that we abuse notation again and denote $C^r$ norms, smoothing operators and average of tuples the same as in the case of smooth maps on torus. With this definitions, the estimates (1.2) still hold in the case of tuples of maps on the torus. We define $C^r$ norms for matices of smooth maps similarly.

We set $\mathcal{F} = (f_1, \ldots, f_m)$, $\mathcal{R}_\alpha = (R_{\alpha_1}, \ldots, R_{\alpha_m})$ and the error term $\tilde{\mathcal{F}} = \mathcal{F} - \mathcal{R}_\alpha = (\tilde{f}_1, \ldots, \tilde{f}_m)$. We define an operator $\mathcal{H} \colon C^\infty(\mathbb{T}^d) \to C^\infty(\mathbb{T}^d)^m$ by $\mathcal{H}(h) = (h, \ldots, h)$. Usually, if it is clear from the context, we just write $\mathcal{H}$ instead of $\mathcal{H}(h)$ and $\tilde{\mathcal{H}}$ instead of $\mathcal{H}(h - id)$. We write $\mathcal{H}^{-1}$ to denote $\mathcal{H}(h^{-1})$, provided $h^{-1}$ exists. We also use the same letter $\mathcal{H}$ for the operator, defined in the same way, on the space of operators acting on $C^\infty(\mathbb{T}^d)$, i.e., if $\mathcal{S} \colon C^\infty(\mathbb{T}^d) \to C^\infty(\mathbb{T}^d)$ is an operator, then $\mathcal{H}(\mathcal{S}) = (\mathcal{S}, \ldots, \mathcal{S})$ stands for the operator on $C^\infty(\mathbb{T}^d)^m$ acting componentwise.

Here we work with lifts of maps and of their tuples and by abuse of notation, we denote them with the same letters. Since a lift of the rotation $R_\beta$ on $\mathbb{T}^d$ is the rotation $R_\beta$ on $\mathbb{R}^d$, for a lift of $\mathcal{F} = \mathcal{R}_\alpha + \tilde{\mathcal{F}}$ we have $\tilde{\mathcal{F}} \in C^\infty(\mathbb{T}^d, \mathbb{R}^d)^m$, where $C^\infty(\mathbb{T}^d, \mathbb{R}^d)$ stands for the space of of smooth $\mathbb{Z}^d$ periodic maps from $\mathbb{R}^d$ to $\mathbb{R}^d$. Everything defined and stated in the previous two paragraphs also apply to the space $C^\infty(\mathbb{T}^d, \mathbb{R}^d)$.

## 3.2 Method of the proof

The goal is to construct a smooth diffeomorphism $h \colon \mathbb{T}^d \to \mathbb{T}^d$, such that $h \circ f_j = R_{\alpha_j} \circ h$, for $1 \leq j \leq m$. We rewrite this in the following, at least notationally simpler, form.

$$\mathcal{F} \circ \mathcal{H} = \mathcal{H} \circ \mathcal{R}_\alpha. \tag{3.1}$$

For this reason, we define a nonlinear operator $\mathcal{N} \colon C^\infty(\mathbb{T}^d)^m \to C^\infty(\mathbb{T}^d)^m$, $\mathcal{N} = (\mathcal{N}_1, \ldots, \mathcal{N}_m)$, where each $\mathcal{N}_j \colon C^\infty(\mathbb{T}^d) \to C^\infty(\mathbb{T}^d)$ is given by $\mathcal{N}_j(h) = f_j \circ h - h \circ R_{\alpha_j}$. Since $\mathcal{F} = \mathcal{R}_\alpha + \tilde{\mathcal{F}}$, solving (3.1) is equivalent to solving the equation

$$\mathcal{DN}(\tilde{\mathcal{H}}) = \tilde{\mathcal{H}} \circ \mathcal{R}_\alpha - \tilde{\mathcal{H}} = \tilde{\mathcal{F}} \circ \mathcal{H}. \tag{3.2}$$

$\mathcal{DN} = (\mathcal{DN}_1, \ldots, \mathcal{DN}_m)$, where $\mathcal{DN}_j \colon C^\infty(\mathbb{T}^d) \to C^\infty(\mathbb{T}^d)$ is the linearization, around the identity map, of the nonlinear operator $\mathcal{N}_j$. Hence, $\mathcal{DN}_j(\varphi) = \varphi \circ R_{\alpha_j} - \varphi$, for $1 \leq j \leq m$.

The proof is based on KAM (Kolmogorov-Arnold-Moser) iterative scheme. The solution of the nonlinear equation (3.2) is the limit of successive approximations obtained by solving approximately the corresponding linear equation

$$\mathcal{DN}(\tilde{\mathcal{H}}) = \tilde{\mathcal{F}}. \qquad (3.3)$$

## 3.3 Solving the linearized equation approximately

The necessary commutativity conditions are given by

$$f_j \circ f_l = f_l \circ f_j, \text{ for } 1 \leq j, l \leq m. \qquad (3.4)$$

Since $\mathcal{F} = \mathcal{R}_\alpha + \tilde{\mathcal{F}}$, condition (3.4) transforms into

$$\tilde{f}_j \circ f_l - \tilde{f}_j = \tilde{f}_l \circ f_j - \tilde{f}_l, \text{ for } 1 \leq j, l \leq m. \qquad (3.5)$$

We need to find an approximate solution $\tilde{\mathcal{H}}$ of the linear equation

$$\mathcal{DN}(\tilde{\mathcal{H}}) = \tilde{\mathcal{F}},$$

where $\tilde{\mathcal{H}}$ is a tuple of smooth maps with zero average. For this purpose we introduce spaces $C_0^\infty(\mathbb{T}^d, \mathbb{R}^d)$, $C_0^\infty(\mathbb{T}^d, \mathbb{R}^d)^m$ and $C_0^\infty(\mathbb{T}^d, \mathbb{R}^d)^{m \times m}$, where $C_0^\infty(\mathbb{T}^d, \mathbb{R}^d)$ is the space of maps from $C^\infty(\mathbb{T}^d, \mathbb{R}^d)$, whose averages are zero. These spaces inherit $C^r$ norms defined above. We also introduce a linear operator $\mathcal{DM} \colon C_0^\infty(\mathbb{T}^d, \mathbb{R}^d)^m \to C_0^\infty(\mathbb{T}^d, \mathbb{R}^d)^{m \times m}$ as

$$\mathcal{DM}(\mathcal{A})_{jl} = \mathcal{DN}_j a_l - \mathcal{DN}_l a_j. \qquad (3.6)$$

Since all $\mathcal{DN}_j$ commute with each other, one has immediately $\mathcal{DM} \circ \mathcal{DN} = 0$, meaning that the necessary condition, for solving the equation (3.3), is $\mathcal{DM}(\tilde{\mathcal{F}}) = 0$. In next lemma we show that this condition is not only necessary but also sufficient for solving (3.3), provided (2.2) is satisfied.

**Lemma 3.1** *If $\mathcal{DM}(\tilde{\mathcal{F}}) = 0$ and $\alpha_1, \ldots, \alpha_m$ are simultaneously Diophantine of type $(\gamma, \tau)$, where $\tau > d$, then there exists $\tilde{\mathcal{H}} \in C_0^\infty(\mathbb{T}^d, \mathbb{R}^d)^m$ such that $\mathcal{DN}(\tilde{\mathcal{H}}) = \tilde{\mathcal{F}}$.*

**Proof.** The condition $\mathcal{DM}(\tilde{\mathcal{F}}) = 0$ splits into finitely many conditions

$$\tilde{f}_l \circ R_{\alpha_j} - \tilde{f}_l = \tilde{f}_j \circ R_{\alpha_l} - \tilde{f}_j, \text{ for } 1 \leq j, l \leq m. \qquad (3.7)$$

After passing to Fourier series of components, we obtain

$$\frac{\hat{\tilde{f}}_{ls}(k)}{e^{2\pi i \langle k, \alpha_l \rangle} - 1} = \frac{\hat{\tilde{f}}_{js}(k)}{e^{2\pi i \langle k, \alpha_j \rangle} - 1}, \qquad (3.8)$$

for every $k \in \mathbb{Z}^d \setminus \{0\}$, every $s \in \{1, \ldots d\}$ and $1 \leq j, l \leq m$. Set

$$\tilde{h}_j(x) = \sum_{k \in \mathbb{Z}^d \setminus \{0\}} \hat{\tilde{h}}_j(k) e^{2\pi i \langle k, x \rangle}, \qquad (3.9)$$

where $\hat{\tilde{h}}_j(k)$ is given by (3.8). Since (2.2) and (3.8) hold, $\tilde{h} \in C_0^\infty(\mathbb{T}^d, \mathbb{R}^d)$, where $\tilde{h} = (\tilde{h}, \ldots, \tilde{h}_d)$. A straightforward computation gives $\mathcal{DN}(\tilde{\mathcal{H}}) = \tilde{\mathcal{F}}$.

**Remark 3.1** *Note that the previous lemma, together with $\mathcal{DM} \circ \mathcal{DN} = 0$, implies that the sequence*

$$0 \longrightarrow \mathrm{C}_0^\infty(\mathbb{T}^d, \mathbb{R}^d) \xrightarrow{\mathcal{DN}} \mathrm{C}_0^\infty(\mathbb{T}^d, \mathbb{R}^d)^m \xrightarrow{\mathcal{DM}} \mathrm{C}_0^\infty(\mathbb{T}^d, \mathbb{R}^d)^{m \times m} \longrightarrow 0 \quad (3.10)$$

*is exact at $\mathrm{C}_0^\infty(\mathbb{T}^d, \mathbb{R}^d)^m$. Unfortunately, there is no reason why $\tilde{\mathcal{F}}$ should satisfy the condition $\mathcal{DM}(\tilde{\mathcal{F}}) = 0$. This means, that in general one cannot solve the linearized equation (3.3) exactly. However, here the commutativity condition enters the picture. Namely, the condition (3.5) ensures that $\tilde{\mathcal{F}}$ almost satisfy the condition $\mathcal{DM}(\tilde{\mathcal{F}}) = 0$, i.e., it ensures that one can solve the linearized equation (3.3) approximately, and that the approximate solution is good enough.*

We now define two more linear operators, which help us in constructing an approximate solution of the linearized equation (3.3). They are

$$\mathcal{DN}^* \colon \mathrm{C}_0^\infty(\mathbb{T}^d, \mathbb{R}^d)^m \to \mathrm{C}_0^\infty(\mathbb{T}^d, \mathbb{R}^d), \; \mathcal{DM}^* \colon \mathrm{C}_0^\infty(\mathbb{T}^d, \mathbb{R}^d)^{m \times m} \to \mathrm{C}_0^\infty(\mathbb{T}^d, \mathbb{R}^d)^m$$

$$\mathcal{DN}^*(\mathcal{A}) = \sum_{j=1}^m \mathcal{DN}_j^* a_j, \; (\mathcal{DM}^*(\mathcal{B}))_j = \sum_{l=1}^m \mathcal{DN}_l^* b_{lj} \quad (3.11)$$

where $\mathcal{DN}_j^* \colon \mathrm{C}^\infty(\mathbb{T}^d) \to \mathrm{C}^\infty(\mathbb{T}^d)$ is defined by

$$\mathcal{DN}_j^* \varphi = \varphi \circ R_{-\alpha_j} - \varphi. \quad (3.12)$$

Since all $\mathcal{DN}_j$ and $\mathcal{DN}_l^*$ commute, we also have $\mathcal{DN}^* \circ \mathcal{DM}^* = 0$ and for the map

$$\mathcal{DN} \circ \mathcal{DN}^* + \mathcal{DM}^* \circ \mathcal{DM} \colon \mathrm{C}_0^\infty(\mathbb{T}^d, \mathbb{R}^d)^m \to \mathrm{C}_0^\infty(\mathbb{T}^d, \mathbb{R}^d)^m,$$

the following identity holds.

$$\mathcal{DN} \circ \mathcal{DN}^* + \mathcal{DM}^* \circ \mathcal{DM} = \mathcal{H}(\mathcal{DN}^* \circ \mathcal{DN}) \quad (3.13)$$

Anticipating the previous identity, in next lemma we define a linear operator whose inverse we use in the definition of an approximate solution of the linearized equation (3.3), provided it exists. We also give estimates for the inverse operator.

**Lemma 3.2** *The operator*

$$\mathcal{P} = \mathcal{H}(\mathcal{DN}^* \circ \mathcal{DN}) \colon \mathrm{C}_0^\infty(\mathbb{T}^d, \mathbb{R}^d)^m \to \mathrm{C}_0^\infty(\mathbb{T}^d, \mathbb{R}^d)^m$$

*is a bijection if and only if $\alpha_1, \ldots, \alpha_m$ are simultaneously Diophantine of type $(\gamma, \tau)$, where $\tau > d$. In that case, we have the following estimate*

$$\left\| \mathcal{P}^{-1}(\mathcal{A}) \right\|_r \leq C_r \left\| \mathcal{A} \right\|_{r+\sigma}, \quad (3.14)$$

*for every $r \geq 0$, where $\mathcal{A} = (a_1, \ldots, a_m)$ and $\sigma = 2(\tau + 1) + d$.*

Proof. It is enough to prove lemma for the operator

$$\mathcal{DN}^* \circ \mathcal{DN} \colon \mathrm{C}_0^\infty(\mathbb{T}^d, \mathbb{R}^d) \to \mathrm{C}_0^\infty(\mathbb{T}^d, \mathbb{R}^d).$$

Let $g = (g_1, \ldots, g_d) \in \mathrm{C}_0^\infty(\mathbb{T}^d, \mathbb{R}^d)$. For $1 \leq j \leq d$, it holds

$$((\mathcal{DN}^* \circ \mathcal{DN})^{-1} g)_j(x) = \sum_{k \in \mathbb{Z}^d \setminus \{0\}} p_k^{-1} \hat{g}_j(k) e^{2\pi i \langle k, x \rangle}, \quad (3.15)$$

where

$$p_k = 2 \sum_{s=1}^{m} (1 - \cos(2\pi \langle k, \alpha_s \rangle)) = 4 \sum_{s=1}^{m} \sin^2(\pi \langle k, \alpha_s \rangle). \quad (3.16)$$

Since (2.2) holds, we have $p_k > 0$ for all $k \in \mathbb{Z}^d \setminus \{0\}$, and it holds

$$p_k^{-1} \leq \sum_{s=1}^{m} \frac{1}{|\langle k, \alpha_s \rangle|^2} \leq C|k|^{2\tau}. \quad (3.17)$$

Since $\hat{g}_j(k) \leq (2\pi |k|)^{-r} \|g_j\|_r$, using (3.17) we obtain

$$\begin{aligned}
\left\| ((\mathcal{DN}^* \circ \mathcal{DN})^{-1} g)_j \right\|_r &\leq \sum_{k \in \mathbb{Z}^d \setminus \{0\}} p_k^{-1} (2\pi |k|)^r \|g_j\|_r \\
&\leq C_r \sum_{k \in \mathbb{Z}^d \setminus \{0\}} |k|^{2\tau - \sigma} \|g_j\|_{r+\sigma} \\
&\leq C_r \|g_j\|_{r+\sigma},
\end{aligned} \quad (3.18)$$

which then implies (3.14).

**Remark 3.2** *Note that $\mathcal{P}^{-1}$ is not bounded as an operator on the same space, but rather is tamely bounded as an operator on a Frechét space. Thus, we have some loss of regularity. Fortunately, the loss is fixed, and depends only on constants appearing in the simultaneously Diophantine condition (2.2), and dimension of the torus. To overcome this loss of regularity, we use the smoothing operators $T_N$ defined at the end of Section 1, for an appropriate choice of $N \geq 0$.*

We define an approximate solution of the linearized equation (3.3) by

$$\tilde{h} = \mathcal{DN}^* \circ \mathcal{P}^{-1} \circ T_N(\tilde{\mathcal{F}}). \quad (3.19)$$

We note that $\mathcal{P}^{-1}$ is not defined on the whole $C^\infty(\mathbb{T}^d, \mathbb{R}^d)^m$, but rather on $C_0^\infty(\mathbb{T}^d, \mathbb{R}^d)^m$. To handle this we simply extend it linearly by zero.

### 3.4 Estimates of the new error

Now that we have an approximate solution of (3.3), we define $h = id + \tilde{h}$. Assuming that $\left\| \tilde{h} \right\|_1 < \eta$ (where $0 < \eta < 1$ is such that that $h^{-1}$ exists), we also define $\mathcal{F}^{(1)} = \mathcal{H}^{-1} \circ \mathcal{F} \circ \mathcal{H}$. The goal is to show that the new error $\tilde{\mathcal{F}}^{(1)} = \mathcal{F}^{(1)} - \mathcal{R}_\alpha$ is quadratically small with respect to the old error $\tilde{\mathcal{F}}$. Since $\mathcal{H} \circ \mathcal{F}^{(1)} = \mathcal{F} \circ \mathcal{H}$, we have $\mathcal{F}^{(1)} + \tilde{\mathcal{H}} \circ \mathcal{F}^{(1)} = \mathcal{R}_\alpha + \tilde{\mathcal{H}} + \tilde{\mathcal{F}} \circ \mathcal{H}$ and thus

$$\begin{aligned}
\tilde{\mathcal{F}}^{(1)} &= \tilde{\mathcal{H}} - \tilde{\mathcal{H}} \circ \mathcal{F}^{(1)} + \tilde{\mathcal{F}} \circ \mathcal{H} \\
&= -\mathcal{DN}(\tilde{\mathcal{H}}) + \tilde{\mathcal{H}} \circ \mathcal{R}_\alpha - \tilde{\mathcal{H}} \circ \mathcal{F}^{(1)} + \tilde{\mathcal{F}} \circ \mathcal{H} \\
&= -\mathcal{DN}(\tilde{\mathcal{H}}) + \tilde{\mathcal{F}} + \mathcal{E},
\end{aligned} \quad (3.20)$$

where $\mathcal{E} = \tilde{\mathcal{F}} \circ \mathcal{H} - \tilde{\mathcal{F}} + \tilde{\mathcal{H}} \circ \mathcal{R}_\alpha - \tilde{\mathcal{H}} \circ \mathcal{F}^{(1)}$.

We first show that $\mathcal{E}$ is quadratically small with respect to the old error $\tilde{\mathcal{F}}$. For $0 \leq r \leq s$, using (3.14), (3.19) and (1.2), we obtain the estimate

$$\left\|\tilde{\mathcal{H}}\right\|_r \leq C \left\|\mathcal{P}^{-1} \circ T_N(\tilde{\mathcal{F}})\right\|_r \leq C_r \left\|T_N(\tilde{\mathcal{F}})\right\|_{r+\sigma} \leq C_{r,s} N^{r-s+\sigma+\frac{d}{2}} \left\|\tilde{\mathcal{F}}\right\|_s. \tag{3.21}$$

Using this for $r = s = 0$, the first part in $\mathcal{E}$ can be estimated in the following way

$$\left\|\tilde{\mathcal{F}} \circ \mathcal{H} - \tilde{\mathcal{F}}\right\|_0 \leq \left\|\tilde{\mathcal{F}}\right\|_1 \left\|\tilde{\mathcal{H}}\right\|_0 \leq C N^{\sigma+\frac{d}{2}} \left\|\tilde{\mathcal{F}}\right\|_1 \left\|\tilde{\mathcal{F}}\right\|_0,$$

and the second part in $\mathcal{E}$ can be dominated by $\frac{1}{4}\left\|\tilde{\mathcal{F}}^{(1)}\right\|_0$, if we assume that $\left\|\tilde{\mathcal{H}}\right\|_1 < \frac{\eta}{4}$, since $\left\|\tilde{\mathcal{H}} \circ \mathcal{R}_\alpha - \tilde{\mathcal{H}} \circ \mathcal{F}^{(1)}\right\|_0 \leq \left\|\tilde{\mathcal{H}}\right\|_1 \left\|\tilde{\mathcal{F}}^{(1)}\right\|_0$. From here, we conclude that

$$\|\mathcal{E}\|_0 \leq C N^{\sigma+\frac{d}{2}} \left\|\tilde{\mathcal{F}}\right\|_1 \left\|\tilde{\mathcal{F}}\right\|_0.$$

Now we estimate the expression $-\mathcal{DN}(\tilde{\mathcal{H}}) + \tilde{\mathcal{F}}$ in (3.20). Using (3.19) and (3.13), we obtain

$$\begin{aligned}-\mathcal{DN}(\tilde{\mathcal{H}}) + \tilde{\mathcal{F}} &= -\mathcal{DN} \circ \mathcal{DN}^* \circ \mathcal{P}^{-1} \circ T_N(\tilde{\mathcal{F}}) \\ &= \mathcal{DM}^* \circ \mathcal{DM} \circ \mathcal{P}^{-1} \circ T_N(\tilde{\mathcal{F}}) - \mathcal{P} \circ \mathcal{P}^{-1} \circ T_N(\tilde{\mathcal{F}}) + \tilde{\mathcal{F}}.\end{aligned} \tag{3.22}$$

Since we extended linearly $\mathcal{P}^{-1}$ by zero, we have

$$\begin{aligned}\mathcal{P} \circ \mathcal{P}^{-1} \circ T_N(\tilde{\mathcal{F}}) &= \mathcal{P} \circ \mathcal{P}^{-1}(T_N(\tilde{\mathcal{F}}) - Av(\tilde{\mathcal{F}})) \\ &= T_N(\tilde{\mathcal{F}}) - Av(\tilde{\mathcal{F}})\end{aligned} \tag{3.23}$$

Substituting this into (3.20), we get

$$-\mathcal{DN}(\tilde{\mathcal{H}}) + \tilde{\mathcal{F}} = \mathcal{DM}^* \circ \mathcal{DM} \circ \mathcal{P}^{-1} \circ T_N(\tilde{\mathcal{F}}) + (I - T_N)(\tilde{\mathcal{F}}) + Av(\tilde{\mathcal{F}}). \tag{3.24}$$

Using the definition of $\mathcal{DM}$, commutativity condition (3.5) and the mean value theorem, we obtain the following estimate.

$$\begin{aligned}\left\|\mathcal{DM}(\tilde{\mathcal{F}})_{jl}\right\|_0 &= \left\|\mathcal{DN}_j(\tilde{f}_l) - \mathcal{DN}_l(\tilde{f}_j)\right\|_0 \\ &= \left\|\tilde{f}_l \circ R_{\alpha_j} - \tilde{f}_l \circ f_j + \tilde{f}_j \circ R_{\alpha_l} - \tilde{f}_j \circ f_l\right\|_0 \\ &\leq \left\|\tilde{f}_j\right\|_1 \left\|\tilde{f}_l\right\|_0 + \left\|\tilde{f}_l\right\|_1 \left\|\tilde{f}_j\right\|_0.\end{aligned} \tag{3.25}$$

Hence, $\left\|\mathcal{DM}(\tilde{\mathcal{F}})\right\|_0 \leq C \left\|\tilde{\mathcal{F}}\right\|_1 \left\|\tilde{\mathcal{F}}\right\|_0$. Using this and (3.14), we arrive at

$$\left\|\mathcal{DM}^* \circ \mathcal{DM} \circ \mathcal{P}^{-1} \circ T_N(\tilde{\mathcal{F}})\right\|_0 \leq C \left\|\mathcal{DM} \circ \mathcal{P}^{-1} \circ T_N(\tilde{\mathcal{F}})\right\|_0$$
$$\leq C \left\|\mathcal{P}^{-1} \circ T_N \circ \mathcal{DM}(\tilde{\mathcal{F}})\right\|_0$$
$$\leq C \left\|T_N(\mathcal{DM}\tilde{\mathcal{F}})\right\|_\sigma \quad (3.26)$$
$$\leq C N^{\sigma+\frac{d}{2}} \left\|\mathcal{DM}(\tilde{\mathcal{F}})\right\|_0$$
$$\leq C N^{\sigma+\frac{d}{2}} \left\|\tilde{\mathcal{F}}\right\|_1 \left\|\tilde{\mathcal{F}}\right\|_0.$$

Combining the estimates we have learned so far, (3.20) and (3.24), we see that $\tilde{\mathcal{F}}^{(1)}$ is quadratically small with respect to $\tilde{\mathcal{F}}$, aside from the truncation error and the average of $\tilde{\mathcal{F}}$, i.e.,

$$\left\|\tilde{\mathcal{F}}^{(1)}\right\|_0 \leq C N^{\sigma+\frac{d}{2}} \left\|\tilde{\mathcal{F}}\right\|_1 \left\|\tilde{\mathcal{F}}\right\|_0 + \left\|(I - T_N)(\tilde{\mathcal{F}})\right\|_0 + \left\|Av(\tilde{\mathcal{F}})\right\|. \quad (3.27)$$

In [4] is shown that $\left\|Av(\tilde{\mathcal{F}})\right\|$ also obeys a good estimate (Lemma 4.1 and Corollary 4.3 in [4]), with respect to $\tilde{\mathcal{F}}$, provided $\alpha_j \in \text{conv}(\rho(f_j))$, for $1 \leq j \leq m$. We state this in a lemma below.

**Lemma 3.3** *If $\alpha_j \in \text{conv}(\rho(f_j))$, for $1 \leq j \leq m$, are simultaneously Diophantine vectors of type $(C, \tau)$, then for every $r \geq 0$ there exists a constant $C_r > 0$ such that $\left\|Av(\tilde{\mathcal{F}})\right\|$ is bounded above by*

$$C_r \left( N^{2(\tau+1)+d} \left\|\tilde{\mathcal{F}}\right\|_0^2 + N^{-r+\tau+\frac{3d}{2}} \left\|\tilde{\mathcal{F}}\right\|_0 \left\|\tilde{\mathcal{F}}\right\|_r + N^{-r+d} \left\|\tilde{\mathcal{F}}\right\|_r \right). \quad (3.28)$$

**Remark 3.3** Note that the term $\left\|(I - T_N)(\tilde{\mathcal{F}})\right\|_0$ in (3.27), is absorbed by $\left\|Av(\tilde{\mathcal{F}})\right\|$, since for any $s \geq 0$, from (1.2) we obtain

$$\left\|(I - T_N)(\tilde{\mathcal{F}})\right\|_0 \leq C_s N^{-s+d} \left\|\tilde{\mathcal{F}}\right\|_s.$$

Estimates (3.27) and (3.28) put together, are enough to ensure the existence of a solution of the nonlinear equation (3.2), provided $\left\|\tilde{\mathcal{F}}\right\|_l$ is small enough, for some large $l \in \mathbb{N}$.

We also need an estimate on $\left\|\tilde{\mathcal{F}}^{(1)}\right\|_l$, for large enough $l \in \mathbb{N}$. This estimate, for example, can be easily derived from estimates in [5] (Appendix II) and (3.21), for $r = s = l$. Namely, for every $l \geq 1$ there exists a constant $C_l > 0$ such that

$$\left\|\tilde{\mathcal{F}}^{(1)}\right\|_l \leq C_l N^{\sigma+\frac{d}{2}} \left(1 + \left\|\tilde{\mathcal{F}}\right\|_l\right). \quad (3.29)$$

## 3.5 Iterative procedure and its convergence

We first set $\tilde{\mathcal{F}}^{(0)} = \tilde{\mathcal{F}}$, $\mathcal{F}^{(0)} = \mathcal{R}_\alpha + \tilde{\mathcal{F}}^{(0)}$, $\mathcal{H}^{(0)} = \mathcal{H}(id)$ and $\tilde{\mathcal{H}}^{(0)} = 0$. Then we construct inductively the sequence $\tilde{\mathcal{F}}^{(n)}$, for every $n \geq 1$, in the following way. We choose an appropriate positive integer $N_n$, so that, after solving the linearized equation (3.3), we obtain new $\tilde{\mathcal{H}}^{(n)}$. We define

$$\mathcal{H}^{(n)} = \mathcal{H}^{(0)} + \tilde{\mathcal{H}}^{(n)},$$
$$\mathcal{F}^{(n+1)} = (\mathcal{H}^{(n)})^{-1} \circ \mathcal{F}^{(n)} \circ \mathcal{H}^{(n)}, \qquad (3.30)$$
$$\tilde{\mathcal{F}}^{(n+1)} = \mathcal{F}^{(n+1)} - \mathcal{R}_\alpha.$$

Hence, $\mathcal{F}^{(n+1)} = H_n^{-1} \circ \mathcal{F} \circ H_n$, where $H_n = \mathcal{H}^{(0)} \circ \ldots \circ \mathcal{H}^{(n)}$. In order to show the existence of a solution for the nonlinear equation (3.2), we need to show the convergence of this iterative process. For this reason, we set $\epsilon_n = \epsilon^{(\kappa^n)}$, $\kappa = \frac{4}{3}$, $N_n = \epsilon_n^{-\frac{1}{9(\sigma+1+d/2)}}$ and we fix $l = 2(\sigma + 1 + d)$. We show by induction, that for every $n \in \mathbb{N}$ the following inequalities hold.

$$\begin{aligned}\left\|\tilde{\mathcal{F}}^{(n)}\right\|_0 &< \epsilon_n \\ \left\|\tilde{\mathcal{F}}^{(n)}\right\|_l &< \epsilon_n^{-1} \\ \left\|\tilde{\mathcal{H}}^{(n)}\right\|_1 &< \epsilon_n^{1/2}.\end{aligned} \qquad (3.31)$$

Using (3.29), we obtain

$$\left\|\tilde{\mathcal{F}}^{(n+1)}\right\|_l \leq C_l N_n^{\sigma+\frac{d}{2}} \left(1 + \left\|\tilde{\mathcal{F}}^{(n)}\right\|_l\right) \leq C_l N_n^{\sigma+\frac{d}{2}} \left(1 + \epsilon_n^{-1}\right)$$
$$\leq 2C_l N_n^{\sigma+\frac{d}{2}} \epsilon_n^{-1} \leq C_l \epsilon_n^{-\frac{\sigma+d/2}{9(\sigma+1+d/2)}} \epsilon_n^{-1}$$
$$< \epsilon_n^{-1/3-1} = \epsilon_n^{-\kappa} = (\epsilon_{n+1})^{-1}.$$

From interpolation estimates for $C^r$ norms, we have

$$\left\|\tilde{\mathcal{F}}^{(n)}\right\|_1 \leq C_r \left\|\tilde{\mathcal{F}}^{(n)}\right\|_0^{1-\frac{1}{r}} \left\|\tilde{\mathcal{F}}^{(n)}\right\|_r^{\frac{1}{r}},$$

for any $r \geq 1$, and hence,

$$CN_n^{\sigma+\frac{d}{2}} \left\|\tilde{\mathcal{F}}^{(n)}\right\|_1 \left\|\tilde{\mathcal{F}}^{(n)}\right\|_0 \leq C\epsilon_n^{-\frac{\sigma+d/2}{9(\sigma+1+d/2)}} \epsilon_n^{1-\frac{1}{r}} \epsilon_n^{-\frac{1}{r}} \epsilon_n$$
$$= C\epsilon_n^{-\frac{\sigma+d/2}{9(\sigma+1+d/2)} + 2(1-\frac{1}{r})} \qquad (3.32)$$
$$< \epsilon_n^{4/3} = \epsilon_{n+1},$$

Using Lemma 3.3, we estimate each part separately in (3.28). From interpolation inequalites, we get

$$N_n^{-r+d} \left\|\tilde{\mathcal{F}}^{(n)}\right\|_r \leq C_r N_n^{-r+d} \left\|\tilde{\mathcal{F}}^{(n)}\right\|_0^{1-\frac{r}{l}} \left\|\tilde{\mathcal{F}}^{(n)}\right\|_l^{\frac{r}{l}}$$
$$< C_r \epsilon_n^{\frac{r-d}{9(\sigma+1+d/2)}} \epsilon_n^{1-\frac{r}{l}} \epsilon_n^{-\frac{r}{l}}$$
$$= C\epsilon_n^{\frac{r-d}{9(\sigma+1+d/2)} + 1 - \frac{2r}{l}} \qquad (3.33)$$
$$< \epsilon_n^{4/3} = \epsilon_{n+1},$$

for every $\frac{3\sigma+3+5d/2}{2} < r < l$.

Again, from interpolation estimates, we get

$$\begin{aligned}
N_n^{-r+\tau+\frac{3d}{2}} \left\|\tilde{\mathcal{F}}^{(n)}\right\|_0 \left\|\tilde{\mathcal{F}}^{(n)}\right\|_r &\leq C_r N_n^{-r+\tau+\frac{3d}{2}} \left\|\tilde{\mathcal{F}}^{(n)}\right\|_0 \left\|\tilde{\mathcal{F}}^{(n)}\right\|_0^{1-\frac{r}{l}} \left\|\tilde{\mathcal{F}}^{(n)}\right\|_l^{\frac{r}{l}} \\
&\leq C_r \epsilon_n^{-\frac{-r+3d/2+\sigma/2-1}{9(\sigma+1+d/2)}} \epsilon_n^{2-\frac{r}{l}} \epsilon_n^{-\frac{r}{l}} \\
&= C_r \epsilon_n^{-\frac{-r+3d/2+\sigma/2-1}{9(\sigma+1+d/2)}+2(1-\frac{r}{l})} \\
&< \epsilon_n^{4/3} = \epsilon_{n+1},
\end{aligned} \tag{3.34}$$

for every $\frac{l}{l+9} \cdot \frac{\sigma+3d-14}{2} < r < l$.

And finally,

$$\begin{aligned}
N_n^\sigma \left\|\tilde{\mathcal{F}}^{(n)}\right\|_0^2 &\leq C \epsilon_n^{-\frac{\sigma}{9(\sigma+1+d/2)}} \epsilon_n^2 \\
&< C \epsilon_n^{2-\frac{1}{9}} \\
&< \epsilon_n^{4/3} = \epsilon_{n+1}.
\end{aligned} \tag{3.35}$$

Now, combining the estimates (3.27), (3.32), Remark 3.3, and the estimates (3.33), (3.34) and (3.35), we get

$$\left\|\tilde{\mathcal{F}}^{(n+1)}\right\|_0 < \epsilon_{n+1}.$$

Using (3.21), for $r = 1$ and $s = 0$, we obtain

$$\begin{aligned}
\left\|\tilde{\mathcal{H}}^{(n+1)}\right\|_1 &\leq C N_n^{1+\sigma+\frac{d}{2}} \left\|\tilde{\mathcal{F}}^{(n)}\right\|_0 \leq C \epsilon_n^{-\frac{1+\sigma+d/2}{9(\sigma+1+d/2)}} \epsilon_n \\
&= C \epsilon_n^{1-1/9} < \epsilon_n^{2/3} = (\epsilon_{n+1})^{1/2},
\end{aligned}$$

and hence, for $\left\|\tilde{\mathcal{F}}\right\|_l$ small enough, the iterative process converges in $C^1$ norm to a solution $\tilde{\mathcal{H}}$, for which $\left\|\tilde{\mathcal{H}}\right\|_1 < \frac{\eta}{4}$.

At the end, we show that the iterative process converges to $\tilde{\mathcal{H}}$ in any $C^r$ norm. For any $r \in \mathbb{N}$, from (3.29) we obtain

$$\left\|\tilde{\mathcal{F}}^{(n+1)}\right\|_r \leq C_r N_n^{\sigma+\frac{d}{2}} \left(1 + \left\|\tilde{\mathcal{F}}^{(n)}\right\|_r\right) \leq \epsilon_n^{-1/9} \left(1 + \left\|\tilde{\mathcal{F}}^{(n)}\right\|_r\right).$$

From this, we derive

$$\left\|\tilde{\mathcal{F}}^{(n)}\right\|_r \leq C \prod_{j=1}^{n-1} \epsilon_j^{-1/9} \left(1 + \left\|\tilde{\mathcal{F}}\right\|_r\right) \leq C_r \epsilon_n^{-1},$$

where, in the last inequality, $C_r = C \left(1 + \left\|\tilde{\mathcal{F}}\right\|_r\right)$. Using interpolation estimates once again, we get, for any $j \geq 1$ and $r = 3j$

$$\left\|\tilde{\mathcal{F}}^{(n)}\right\|_j \leq \left\|\tilde{\mathcal{F}}^{(n)}\right\|_0^{2/3} \left\|\tilde{\mathcal{F}}^{(n)}\right\|_r^{1/3} < C_j \epsilon_n^{1/3}$$
$$\left\|\tilde{\mathcal{H}}^{(n)}\right\|_j \leq C_j N_n^{\sigma+\frac{d}{2}} \left\|\tilde{\mathcal{F}}^{(n)}\right\|_j \leq C_j \epsilon_n^\delta,$$

where $\delta = \frac{1}{3(\sigma+1+d/2)} > 0$.